\documentclass{article}
\usepackage{amssymb,latexsym,amsmath,amsthm,mathrsfs}
\usepackage{graphicx}
\newcommand{\comment}[1]{}

\begin{document}
\title{New demonstrations about the resolution of numbers into
squares\footnote{Presented to the St. Petersburg Academy on March 24, 1774.
Originally published as
{\em Novae demonstrationes circa resolutionem numerorum in quadrata},
Nova acta eruditorum (1773),
193--211.
E445 in the Enestr{\"o}m index.
Translated from the Latin by Jordan Bell,
Department of Mathematics, University of Toronto, Toronto, Ontario, Canada.
Email: jordan.bell@utoronto.ca}}
\author{Leonhard Euler}
\date{}
\maketitle

1. Though I have been occupied with this argument intensely and often,
still, the demonstration which I gave previously about the resolution of all
numbers into four or fewer squares had not been altogether satisfactory to me.
Thus I have pursued with all the greater 
ardor the demonstration which 
the Celebrated Mr. Lagrange has recently given
in the first volume of the Berlin
{\em M\'emoires}
of this theorem, which I judge to completely settle the matter,
even if it seems to demand much effort and be very laborious.

2. I believe however that it will be hardly 
unwelcome to readers if I relate briefly and clearly
the particular points on which Lagrange's demonstration
is based. After the Celebrated Author gives the lemma that if
two sums of two squares $pp+qq$ and $rr+ss$ have a common divisor
$\varrho$ which does not divide all the individual squares,
then not
only this divisor $\varrho$ itself, but also both the quotients
$\frac{pp+qq}{\varrho}$ and $\frac{rr+ss}{\varrho}$ will be sums
of two squares, he proceeds to the theorem to be demonstrated,
{\em that if a sum of four squares $P^2+Q^2+R^2+S^2$ is divisible by any
number $A$ which does not divide all the individual squares,
 then the number
$A$ itself will be a sum of four squares}, whose demonstration
is contained in the following reasoning. 

I. Putting the quotient arising from this division $=a$, so that
\[
Aa=P^2+Q^2+R^2+S^2,
\]
if it turns out that the two formulas $P^2+Q^2$ and $R^2+S^2$ have a common
divisor $\varrho$, then since the number $a$ will also contain
it\footnote{Translator: Euler is implicitly taking $A$ prime, and
indeed he only uses this result for $A$ prime. This is false for $A$ composite.
Playing around in Maple I found the
example $P=219, Q=192, R=255, S=402$, and $A=3^2\cdot 13$. $A$ divides
$P^2+Q^2+R^2+S^2$ but does not divide all of $P,Q,R,S$; in fact it divides
none. However, both $P^2+Q^2$ and $R^2+S^2$ are divisible by
$\varrho=3$. Thus here it does not follow that $\varrho$ divides
$a=2\cdot 11^2$.}, one puts
$a=b\varrho$. Thus it becomes
\[
Ab=\frac{P^2+Q^2}{\varrho}+\frac{R^2+S^2}{\varrho};
\]
since by the stated lemma\footnote{Translator: This is referring to the result
stated in \S 2.} these formulas are sums of two squares,
an equation of this kind will be obtained
\[
Ab=pp+qq+rr+ss,
\]
where the formulas $pp+qq$ and $rr+ss$ will no longer have a common factor.

II. Then indeed it is put $pp+qq=t$ and $rr+ss=u$, so that $Ab=t+u$,
which by multiplying by $t$ leads to the equation $Abt=tt+tu$; and because
$tu$ is also a sum of two squares, call it say $xx+yy$, by taking namely
$x=pr+qs$ and $y=ps-qr$ it will become
\[
Abt=tt+xx+yy.
\]

III. Now it is observed that both $x$ and $y$ can thus be expressed
by the numbers $t$ and $b$,  which of course are mutually prime, as
$x=\alpha t+\gamma b$ and $y=\beta t+\delta b$;
where although there are infinitely many ways the letters $\alpha,\beta,\gamma,\delta$ can be taken as either negative or positive, among
which certain values are given so that $\alpha<\frac{1}{2}b$ and
$\beta<\frac{1}{2}b$.

IV. Now substituting these values for $x$ and $y$, this equation will result
\[
Abt=tt(1+\alpha\alpha+\beta\beta)+2bt(\alpha\gamma+\beta\delta)+bb(\gamma\gamma+\delta\delta).
\]
While this expression should be divisible by $b$, however in the first
term $tt$ does not admit this division\footnote{Translator: Since $Ab=t+u$
where $t$ and $u$ are relatively prime.}, so it is necessary that
the formula $1+\alpha\alpha+\beta\beta$ has the factor $b$;
as well, in the same way it is necessary in the last term that the
factor $\gamma\gamma+\delta\delta$ be divisible by $t$.
Let it therefore be put $1+\alpha\alpha+\beta\beta=ba'$, and
because each number $\alpha$ and $\beta$ is less than
$\frac{1}{2}b$, it is clear that $a'<\frac{1}{2}b+\frac{1}{b}$;
then dividing by $b$ it will be
\[
At=a'tt+2t(\alpha\gamma+\beta\delta)+b(\gamma\gamma+\delta\delta).
\]

V. Now let this equation be multiplied by $a'$, so that it becomes
\[
Aa't={a'}^2tt+2a't(\alpha\gamma+\beta\delta)+a'b(\gamma\gamma+\delta\delta),
\]
and in the last term by writing $1+\alpha\alpha+\beta\beta$ in place of
$a'b$ it becomes
\[
Aa't={a'}^2tt+2a't(\alpha\gamma+\beta\delta)+(\alpha\alpha+\beta\beta)(\gamma\gamma+\delta\delta)+\gamma\gamma+\delta\delta.
\]
This is resolved into four squares in the following
\[
Aa't=(a't+\alpha\gamma+\beta\delta)^2+(\beta\delta-\alpha\delta)^2+\gamma^2+\delta^2;
\]
where 
since the sum of the latter two squares $\gamma^2+\delta^2$ is divisible by
the number $t$, it is necessary that the sum of the first
two is also divisible by $t$,
so that here two sums of two squares occur having a common divisor $t$;
whence if they are divided by $t$, both their quotients will likewise
be sums of two squares.

VI. But if we thus put
\[
\frac{(a't+\alpha\gamma+\beta\delta)^2+(\beta\gamma-\alpha\delta)^2}{t}={p'}^2+{q'}^2 \quad \textrm{and} \quad \frac{\gamma^2+\delta^2}{t}={r'}^2+{s'}^2,
\]
we will have
\[
Aa'={p'}^2+{q'}^2+{r'}^2+{s'}^2.
\]
As well, in this formula $Aa'$, if it is compared with the first
$Aa$, the number $a'$ will be much smaller than $a$,
since $b<a$ and $a'<\frac{1}{2}b$.
In a similar way, we can get to a formula $Aa''$, where $a''$ will be much
smaller than $a'$, and thus finally it is necessary that the
formula $A\cdot 1$ be reached, so that now the number $A$ is found to be equal
to a sum of four squares.

3. In order to demonstrate this theorem it is also necessary to show that
for any given prime number a sum of four squares can be exhibited which is divisible by that prime, 
but which are not all divisible by the prime.
And indeed the Celebrated Lagrange
also demonstrates this in a very 
ingenious way, which
however is abstruse and lengthy, so that his efforts 
are not revealed as briefly and clearly as could be desired.
Now therefore, the famous theorem of Bachet or Fermat, {\em that
any number can be resolved into four squares}, is obtained by a perfect
demonstration. 
Since for any prime number a sum of four squares can always be given which
is divisible by it, all prime numbers will be sums of four or fewer squares,
and because it was formerly already demonstrated that the product of two
or more numbers, each of which is a sum of four or fewer squares, can itself
also be divided into four squares, it has now been established
most securely that all numbers whatsoever are the sum of four or fewer squares.

4. Although to be sure it would be a sin to remove anything from the solidity
and rigor of these demonstrations, nevertheless no one will deny that the
foundations and rules of all the reasoning by which these demonstrations are
composed is rather lengthy, and I have not been able to fully remove the
obscurity involved, so that even still clearer and easier demonstrations
could be desired. 

5. After I had carefully studied this new argument, new and rather straight forward
demonstrations of these theorems occurred to me, 
which are useful in this pursuit. The publication of these new demonstrations
surely seems worthwhile; I will relate these here as briefly and clearly
as I am able.
I shall first take up that well known and fully demonstrated theorem that
all divisors of a sum of two relatively prime squares are themselves
equal to a sum of two squares, since this new demonstration commends itself
by its simplicity. Then by following in the same footsteps
the demonstration will be easily extended to four squares.

\begin{center}
{\Large Lemma 1}
\end{center}

6. {\em A product of two sums of two squares is itself a sum of two squares.}

For if that product were $(aa+bb)(\alpha\alpha+\beta\beta)$ and
one takes
\[
A=a\alpha+b\beta \quad \textrm{and} \quad B=a\beta-b\alpha,
\]
it will certainly be
\[
(aa+bb)(\alpha\alpha+\beta\beta)=AA+BB.
\]

\begin{center}
{\Large Theorem 1}
\end{center}

{\em If a number $N$ is a divisor of a sum of two squares $P^2+Q^2$
which are prime to each other, then that number $N$ will itself be a sum
of two squares.}

\begin{center}
{\Large Demonstration}
\end{center}

In order to work out this demonstration with more manageable numbers, I observe
that however large the numbers $P$ and $Q$ are, from
them a sum of two squares $pp+qq$ can always be formed whose roots
$p$ and $q$ do not exceed half of the given number $N$.
For if one puts
\[
P=fN\pm p \quad \textrm{and} \quad Q=gN \pm q,
\]
it is familiar that numbers $p$ and $q$ can be taken as not to exceed
the half $\frac{1}{2}N$. Now since it would then be
\[
PP+QQ=NN(ff+gg)+2N(\pm fp\pm gq)+pp+qq
\]
and this expression would be divisible by $N$, it is evident that this
sum of two squares is also divisible by $N$. With this stated, I will complete
the demonstration by the following steps.

I. Since this formula $pp+qq$ has the divisor $N$, by putting the quotient
$=n$ we will have
\[
Nn=pp+qq,
\]
where therefore $n$ will be less than $\frac{1}{2}N$,
because $p<\frac{1}{2}N$ and $q<\frac{1}{2}N$.

II. Now these numbers $p$ and $q$ can be expressed by the number $n$ such that
\[
p=a+\alpha n \quad \textrm{and} \quad q=b+\beta n,
\] 
where by also admitting negative numbers for $a$ and $b$ it will be
possible to reduce these below $\frac{1}{2}n$, as we observed initially.
Then indeed it will be
\[
N=naa+bb+2n(a\alpha+b\beta)+nn(\alpha\alpha+\beta\beta),
\]
and because in the stated lemma it was $a\alpha+b\beta=A$, it becomes
\[
Nn=aa+bb+2nA+nn(\alpha\alpha+\beta\beta).
\]

III. Therefore the first part $aa+bb$ of this expression necessarily
has a factor $n$, because the other part just by itself admits the divisor
$n$. Therefore let us set
\[
aa+bb=nn',
\]
and because $a<\frac{1}{2}n$ and $b<\frac{1}{2}n$ and then $nn'<\frac{1}{2}nn$,
it will obviously be $n'<\frac{1}{2}n$.
Substituting in this value and dividing by $n$ yields
\[
N=n'+2A+n(\alpha\alpha+\beta\beta).
\]

IV. We multiply this equation by $n'$, and because
$nn'=aa+bb$, by the stated lemma the latter part reduces to 
\[
nn'(\alpha\alpha+\beta\beta)=(aa+bb)(\alpha\alpha+\beta\beta)=AA+BB,
\]
so that we will now have
\[
Nn'=n'n'+2n'A+AA+BB;
\]
this expression is clearly a sum of two squares, namely
\[
Nn'=(n'+A)^2+B^2.
\]

V. Thus as at first it had been that the product $Nn$ was a sum
of two squares and from it we elicited a smaller product $Nn'$ also
equal to a sum of two squares, in this way continually smaller
products can be led to, namely $Nn'', Nn'''$ etc. Thus it is necessary
that finally a minimum product, namely $N\cdot 1$, is reached,
and thus this given number $N$ will also be a sum of two squares.

\begin{center}
{\Large Corollary}
\end{center}

It seems perhaps wonderful that when a number of the type $n'=1$ 
has been reached, the same operations can be applied;\footnote{Translator: My best translation of this paragraph is that it is remarkable how this method of descent stops once we hit $n'=1$, and
that Euler is showing explicitly what happens once we hit $n'=1$.}
this is easily seen by taking $n=1$,
for then one will obtain $p=a+\alpha\cdot 1$ and
$q=b+\beta\cdot 1$, where one clearly takes $a=0$ and $b=0$,
which of course are $<\frac{1}{2}$; then indeed as $aa+bb=0$, it will of course
be $n'=0$ and so this last step spontaneously ends our calculation.

\begin{center}
{\Large Scholion}
\end{center}

It can be demonstrated in the same way that all numbers of
the form $pp+2qq$ or $pp+3qq$ do not admit any other divisors except
those of the same form, if indeed the numbers $p$ and $q$ are prime to
each other. In truth however this calculation cannot be extended to higher forms,
such as $pp+5qq, pp+6qq$, because then
the number $n'$ that is sought is no longer necessarily less
than $n$. Thus let us take up here the demonstrations of the former cases.

\begin{center}
{\Large Lemma 2}
\end{center}

7. {\em A product of two numbers of the form $pp+2qq$ is always a number of this
same form.}

For if such a product is given $(aa+2bb)(\alpha\alpha+2\beta\beta)$
and one takes
\[
A=a\alpha+2b\beta \quad \textrm{and} \quad B=a\beta-b\alpha,
\]
then it will of certainly be
\[
AA+2BB=(aa+2bb)(\alpha\alpha+2\beta\beta).
\]

\begin{center}
{\Large Theorem 2}
\end{center}

{\em If $N$ is a divisor of the number $pp+2qq$, and $p$ and $q$ are
prime to each other, then this number $N$ will also be contained
in such a form.}

\begin{center}
{\Large Demonstration}
\end{center}

Here again the numbers $p$ and $q$ can be kept below half of
the number $N$, and our demonstration will proceed in the following way.

I. Let 
\[
Nn=pp+2qq,
\]
and because $p<\frac{1}{2}N$ and $q<\frac{1}{2}N$, it will
be $n<\frac{3}{4}N$. Now let us put as before
\[
p=a+\alpha n \quad \textrm{and} \quad q=b+\beta n,
\]
where $a$ and $b$ can be taken less than $\frac{1}{2}N$, and then we will
have
\[
Nn=aa+2bb+2n(a\alpha+2b\beta)+nn(\alpha\alpha+2\beta\beta);
\]
by the previous lemma this form is reduced to
\[
Nn=aa+2bb+2nA+nn(\alpha\alpha+\beta\beta).
\]

II. Here therefore the first part $aa+2bb$ will have a factor $n$, whence by putting
\[
aa+2bb=nn'
\]
it will certainly be $n'<\frac{3}{4}n$. Now by substituting this
value and by dividing by $n$ it will become
\[
N=n'+2A+n(\alpha\alpha+2\beta\beta).
\]

III. Let us multiply by $n'$, and then by the previous lemma we will have
\[
nn'(\alpha\alpha+\beta\beta)=(aa+2bb)(\alpha\alpha+2\beta\beta)=AA+2BB,
\]
so that we shall now have
\[
Nn'=n'n'+2n'A+AA+2BB;
\]
this form clearly reduces to 
\[
Nn'=(n'+A)^2+2BB,
\]
and hence likewise a number of the form $pp+2qq$.

IV. Therefore since $n'<n$, in the same way successive products $Nn'',
Nn'''$ etc. can be reached such that the numbers $n,n',n'',n'''$ etc.
continually decrease. Therefore it is at last necessary to reach the
form $N\cdot 1$, so that the number $N$ is itself contained in the same
form $pp+2qq$ too.

\begin{center}
{\Large Lemma 3}
\end{center}

8. {\em A product of two numbers of the form $pp+3qq$ can always be reduced
to the same form.}

For let such a product be $(aa+3bb)(\alpha\alpha+3\beta\beta)$ and take
\[
A=a\alpha+3b\beta \quad \textrm{and} \quad B=a\beta-b\alpha;
\]
one will clearly have
\[
AA+3BB=(aa+3bb)(\alpha\alpha+3\beta\beta).
\]

\begin{center}
{\Large Theorem 3}
\end{center}

{\em If $N$ is a divisor of the number $pp+3qq$, where $p$ and $q$
are numbers which are prime to each other, then this number $N$
will always be able to be reduced to the same form.}

\begin{center}
{\Large Demonstration}
\end{center}

Since again it can be considered $p<\frac{1}{2}N$ and $q<\frac{1}{2}N$,
this form $pp+3qq$ will be less than $N^2$. Then by putting
\[
pp+3qq=Nn
\]
the factor $n$ will be less than $N$,
though in fact this reduction is not necessary for the demonstration;
for it can proceed the same, even if it were $n>N$, as follows.

I. Now by putting
\[
p=a+\alpha n \quad \textrm{and} \quad q=b+\beta n
\]
these numbers $a$ and $b$ can be set less than $\frac{1}{2}n$, at least
not greater;
then it will further be
\[
Nn=aa+3bb+2n(a\alpha+3b\beta)+nn(\alpha\alpha+\beta\beta),
\]
which by the previous lemma is
\[
Nn=aa+3bb+2nA+nn(\alpha\alpha+3\beta\beta).
\]

II. It is therefore necessary that the first part $aa+3bb$ have a factor
$n$; whence by putting
\[
aa+3bb=nn'
\]
this number $n'$ will surely be less than $n$, at least not greater;
then indeed carrying out division by $n$ yields
\[
N=n'+2A+n(\alpha\alpha+3\beta\beta).
\]

III. Now let us multiply by $n'$, and the latter part
\[
nn'(\alpha\alpha+3\beta\beta)=(aa+3bb)(\alpha\alpha+3\beta\beta)
\]
by the previous lemma is $AA+3BB$, and thus we will have
\[
Nn'=n'n'+2n'A+AA+3BB;
\]
this expression clearly reduces to
\[
Nn'=(n'+A)^2+3BB.
\]

IV. Therefore since $Nn'$ is again of the form $pp+3qq$ and $n'<n$,
in the same way continually smaller products $Nn'',Nn'''$ etc. can be
advanced to,
until at last the last $N\cdot 1$ is reached;
and thus it is demonstrated that this number $N$ is  of the form
$pp+3qq$.  

\begin{center}
{\Large Corollary 1}
\end{center}

The basis of this demonstration, like for the preceding,
consists in that for any number $n$, another smaller
$n'$ is reached, which is clear by itself in those cases where $n$ is large
enough. This rule even works in the case where $n=1$; for then
one can take $a=0$ and $b=0$, whence $nn'=0$ which will make it $n'=0$.

Nevertheless clearly a singular case occurs for this theorem when
two ends up being reached in the progression of numbers $n, n', n''$ etc.;
this case merits more attention
because it does not otherwise occur.

\begin{center}
{\Large Corollary 2}
\end{center}

Therefore in the first case let us put $n=2$ and it is clear
that 
in the formula $pp+3qq$ both the numbers $p$ and $q$ must be odd;
for both cannot be assumed to be even, since $p$ and $q$ have
been assumed to be prime to each other. And since here it is
$p=a+2\alpha$ and $q=b+2\beta$, it will become $a=1$ and $b=1$
and hence $aa+3bb=4=nn'$, from which it is clear that $n'$
will be $=2$, so that no further diminution 
can occur here. 

\begin{center}
{\Large Corollary 3}
\end{center}

Here this might become more clear if we reflect that the formula
$pp+3qq$, when both the numbers $p$ and $q$ are odd,
not only cannot be even, but also is divisible by $4$,
thus no oddly even number can be of the form $pp+3qq$.
Therefore whenever, as occurs in these cases, the number $2N$ is contained
in the form $pp+3qq$, then $N$ will always be an even number of which
half, $\frac{1}{2}N$, or a quarter part of $2N$, 
will be contained
in the form $pp+3qq$. For whenever both of the numbers $p$ and $q$
are odd, then too $\frac{pp+3qq}{4}$ is always a number of the
same form, and even in integers, which is not as easy to see.
For by putting $p=2r+1$ and $q=2s+1$,
this formula follows
\[
\frac{pp+3qq}{4}=1+r+rr+3s+3ss,
\]
which can not in general be reduced in integers to a square and a triple
square.\footnote{Translator: I think Euler means that $1+r+rr+3s$
is not necessarily a square, for example $r=1$ and $s=1$ gives
$6+3\cdot 1$, which is not a square by three times a square. But of course
$9=3^2+3\cdot 0$.}
However this resolution can be done in general in the following way.
For I observe that all odd squares are contained in the form
$(4m+1)^2$,
if indeed negative numbers are also admitted for $m$;
if on the one hand $m$ were positive, the squares of the numbers $1,5,9,13$
etc., which are of the form $4i+1$, result;
on the other if $m$ were a negative number, then the squares
of the numbers $3,7,11,15$ etc., which are of the form
$4i-1$, arise. Now let us put
\[
pp=(4r+1)^2 \quad \textrm{and} \quad qq=(4s+1)^2
\]
and it will be
\[
\frac{pp+3qq}{4}=1+2r+4rr+6s+12ss,
\]
which clearly can be put into this form
\[
(1+r+3s)^2+3(r-s)^2.
\]

\begin{center}
{\Large Scholion}
\end{center}

Let us now advance to the demonstration of the stated theorems,\footnote{Translator: In \S 2 of the paper Euler stated 
the results about four squares he is going to prove. Up to now Euler
has been proving results about two squares.} especially
that which is our main object, 
{\em that a sum of four squares admits no other divisors except those which
are a sum of four squares}. 
Like the preceding theorems
it will also be useful to give this lemma. 

\begin{center}
{\Large Lemma 4}
\end{center}

9. {\em A product of two or more numbers, each of which are a sum of
four squares, can always be expressed as a sum of four squares.}

Let such a product be
\[
(aa+bb+cc+dd)(\alpha\alpha+\beta\beta+\gamma\gamma+\delta\delta)
\]
and let us take
\begin{eqnarray*}
A&=&a\alpha+b\beta+c\gamma+d\delta,\\
B&=&a\beta-b\alpha-c\delta+d\gamma,\\
C&=&a\gamma+b\delta-c\alpha-d\beta,\\
D&=&a\delta-b\gamma+c\beta-d\alpha
\end{eqnarray*}
and the sum of the squares of these will be
\[
A^2+B^2+C^2+D^2=(a^2+b^2+c^2+d^2)(\alpha^2+\beta^2+\gamma^2+\delta^2);
\]
for it is clear that all the products of two parts destroy each other,
and all the squares of the Latin letters 
are multiplied by all the squares of the Greek letters.

\begin{center}
{\Large Theorem 4}
\end{center}

{\em If $N$ is a divisor of any sum of four squares, that is of the form
$pp+qq+rr+ss$, each of which indeed are not divisible by $N$, then $N$
will certainly be a sum of four squares.}

\begin{center}
{\Large Demonstration}
\end{center}

It will be of no small help to have noted that these four roots $p,q,r,s$
can be kept below half the given number $N$;
then the demonstration proceeds in the following way.

I. With $n$ denoting the quotient from dividing by this, so that
\[
Nn=pp+qq+rr+ss,
\]
where the letters $p,q,r,s$ may be related thus to $n$
\[
p=a+n\alpha, \quad q=b+n\beta, \quad r=c+n\gamma, \quad s=d+n\delta,
\]
where it is completely obvious that the letters $a,b,c,d$ can be taken so as
not to exceed $\frac{1}{2}n$, since negative values are not excluded here.
And so the formula $aa+bb+cc+dd$ will certainly be less than $nn$.

II. With these values substituted into our equation, it will be
\[
Nn=aa+bb+cc+dd+2n(a\alpha+b\beta+c\gamma+d\delta)+nn(\alpha\alpha+\beta\beta+\gamma\gamma
+\delta\delta),
\]
which by the stated lemma, where we put
\[
A=a\alpha+b\beta+c\gamma+d\delta,
\]
is contracted thus
\[
Nn=aa+bb+cc+dd+2nA+nn(\alpha\alpha+\beta\beta+\gamma\gamma+\delta\delta).
\]
Thus because the first part $aa+bb+cc+dd$ 
should have a factor $n$, let us put
\[
aa+bb+cc+dd=nn'
\]
and it will certainly be $n'<n$, like we just showed. By dividing by $n$ we will
then obtain
\[
N=n'+2A+2n(\alpha\alpha+\beta\beta+\gamma\gamma+\delta\delta).
\]

III. Let us now multiply by $n'$, and because $nn'=aa+bb+cc+dd$, we will have
from the 
preceding lemma
\[
nn'(\alpha\alpha+\beta\beta+\gamma\gamma+\delta\delta)=A^2+B^2+C^2+D^2;
\]
having introduced this form, our equation will be
\[
Nn'=n'n'+2n'A+A^2+B^2+C^2+D^2,
\]
which clearly can be reduced to these four squares
\[
Nn'=(n'+A)^2+B^2+C^2+D^2.
\]

IV. Therefore since $n'<n$, in the same way we can reach continually smaller
forms $Nn'', Nn'''$ etc., until finally we arrive at the form
$N\cdot 1$ and hence the given number $N$ is equated to four squares.

\begin{center}
{\Large Corollary 1}
\end{center}

This calculation is again guilty of a minor exception, namely whenever
it is $n=2$ and all the numbers $p,q,r,s$ are odd;
for then it will happen that $a=1,b=1,c=1$ and $d=1$, so it would also
happen that $n'=2$ and so not less than $n$. 
Truly when the number $2N$ is equal to a sum of four squares,
it is clear from elsewhere that half of it, $N$, is a sum
of four squares, and thus that this exception
should clearly not be considered to disturb anything.

\begin{center}
{\Large Corollary 2}
\end{center}

So that we can see this clearly, let the numbers $p,q,r,s$ be odd and the
number $n$ be even; then, because $Nn=pp+qq+rr+ss$, it will be
\[
\frac{1}{2}Nn=\Big(\frac{p+q}{2}\Big)^2+\Big(\frac{p-q}{2}\Big)^2+
\Big(\frac{r+s}{2}\Big)^2+
\Big(\frac{r-s}{2}\Big)^2,
\]
and these four squares are likewise integers; it will be possible to use
this reduction as long as all the roots of the four squares are odd; for then
the exception mentioned before falls down by itself.

\begin{center}
{\Large Scholion}
\end{center}

By this demonstration this great theorem of Fermat is completed,
since the other part which is still left, namely that given any
number a sum of four squares can be exhibited which is divisible by it,
has been obtained clearly enough by myself for a while already,
and has recently been confirmed by a most exact demonstration by the Celebrated
Lagrange. However so that I can thoroughly complete this argument,
I will adjoin the following very easy demonstration.

\begin{center}
{\Large Theorem 5}
\end{center}

10. {\em Given any prime number $N$, not only four squares but even in fact
three squares can be exhibited in infinitely many ways whose sum is
divisible by this number $N$, but no single one can be divided by it.}

\begin{center}
{\Large Demonstration}
\end{center}

With respect to the number $N$, clearly all numbers are contained
in one of the following forms
\[
\lambda N, \quad \lambda N+1, \quad \lambda N +2, \quad \lambda N+3, \quad
\ldots, \quad \lambda N+N-1,
\]
the number of which is $N$. But disregarding the first form, which contains
the multiples of $N$, it is noted for the remaining, the number of which is
$N-1$, that the squares of the first form $\lambda N+1$ and the last
$\lambda N+N-1$ reduce to the same form $\lambda N+1$; and indeed
the squares of the second form $\lambda N+2$ and the second to last
form $\lambda N+N-2$ to the form $\lambda N+4$; and indeed the third and
third to last can be reduced to $\lambda N+9$, and so on, so that the
squares can themselves be covered by these forms
\[
\lambda N+1,\quad \lambda N+4,\quad \lambda N+9\quad \textrm{etc.},
\]
the number of which is $\frac{1}{2}(N-1)$,
which we shall call forms of the first class and designate thus
\[
\lambda N+a,\quad \lambda N+b,\quad \lambda N+c,\quad \lambda N+d\quad \textrm{etc.},
\]
so that the letters $a,b,c,d$ etc. would denote either the squares
$1,4,9,16$ etc. themselves, or, if they exceed the number $N$, the residues
left from division. Indeed let the other forms, the number of which
will also be $\frac{1}{2}(N-1)$, be designated in this way
\[
\lambda N+\alpha,\quad \lambda N+\beta,\quad \lambda N+\gamma,\quad \lambda N+\delta\quad \textrm{etc.},
\]
which we will call forms of the second class. Let the 
following three properties about these pair of classes be noted, which
indeed can be easily demonstrated.

I. Products of two numbers from the first class are again contained in
the first class, namely the form $\lambda N+ab$ which occurs in the first
class; for if $ab$ were greater than $N$, then the residue resulting from
division by $N$ is to be understood as taken in its place. 

II. Numbers of the first class $a,b,c,d$ etc. multiplied by any number of
the second class $\alpha,\beta,\gamma,\delta$ etc. will end up in the second class.

III. Finally, the product of two numbers from the second class, such as 
$\alpha \beta$, will be transferred to the first class.

With these established I will demonstrate: If three squares could not
be given whose sum were divisible by $N$, then a great absurdity will follow
from this. 
For let, let us concede for the moment the contrary position
that three squares cannot be given whose sum is divisible by $N$;
much less therefore could two such squares be given.
Then it will follow at once that the form $\lambda N-a$, or what reduces
to the same, $\lambda N+(N-a)$, will not appear in the first class; for
given a square of the form $\lambda N-a$, a square of the form $\lambda N+a$
would yield a sum divisible by $N$, contrary to the hypothesis. Therefore,
so that the form $\lambda N-a$ is contained in the latter class
it is necessary that the $\alpha,\beta,\gamma,\delta$ etc. are comprised
 the numbers $-1,-4,-9$ etc.
Let $f$ be any number of the form class, so that a square of the form $\lambda N+f$ are given; if squares of the form $\lambda +1$ 
are added to this, $\lambda N+f+1$ will be a sum of two squares. Now if
a square of the form $\lambda N-f-1$ were given, 
we would obtain a sum of square squares that was divisible by $N$; since
this is false, the form $\lambda N-f-1$ will not be in the first
class and will this be contained in the latter; therefore, since the numbers
$-1$ and $-f-1$ are there, it is necessary that their product $+f+1$ occurs
in the first class. It can be shown in a similar way that these numbers also
must occur in the first class
\[
f+2, \quad f+3, \quad f+4 \quad \textrm{etc.};
\]
whence taking $f=1$, clearly all the forms
\[
\lambda N+1, \quad \lambda N+2, \quad \lambda N+3 \quad \textrm{etc.}
\]
must occur in the first class, and none at all remain for the latter
class. On the other hand, we have seen by the same reasoning that 
the numbers
\[
-1, \quad -f-1, \quad -f-2 \quad \textrm{etc.}
\]
occur in the latter class, and thus clearly all the forms also appear here;
since this is completely absurd, it follows as a falsehood that three squares
cannot be given whose sum is divisible by the given number $N$. Thus there 
are given three, and also four, squares, whose sum will be divisible by $N$.

\begin{center}
{\Large Corollary}
\end{center}

From this theorem in conjunction with the preceding, it follows clearly
that all prime numbers whatsoever 
are sums of four or fewer squares. And since the product of two or more
of these numbers also have this same property, 
it has been most completely shown that {\em all numbers whatsoever are
the sum of four squares or even fewer.}

\begin{center}
{\Large Scholion}
\end{center}

In place of this proposition, the Celebrated Lagrange
gave to the public a theorem holding more widely and supported it by an
ingenious demonstration which however was abstruse and hard to understand,
that could only be understood with the greatest attention.
He showed namely that given any prime number $A$, two squares $pp$ and $qq$
can always be given so that the formula $pp-Bqq-C$ is divisible by the same
prime number $A$, whatever numbers are taken for the letters $B$ and $C$,
providing that they are prime with respect to $A$. 
I will therefore adjoin here the same theorem extended somewhat more widely,
with a much easier and straightforward demonstration.

\begin{center}
{\Large Theorem 6}
\end{center}

{\em 11. Given any prime number $N$, three squares $xx,yy$ and $zz$ prime
to it can always be exhibited so that the formula
\[
\lambda xx+\mu yy+\nu zz
\]
becomes divisible by this prime number $N$, providing these coefficients
$\lambda,\mu$ and $\nu$ are prime to $N$, that is, none of them vanish
and none of them can be made equal either to $N$ itself or any multiple of $N$.}

\begin{center}
{\Large Demonstration}
\end{center}

Let the letters
\[
a, \quad b, \quad c, \quad d \quad \textrm{etc.}
\]
denote all the residues which remain from dividing squares by the given
prime number $N$, which we previously
called the first class,
whose multitude is $\frac{1}{2}(N-1)$; in these namely all the square numbers 
$1,4,9,16$ etc. less than $N$ occur,
together with the residues which are left when the larger
ones are divided by $N$.
Indeed the numbers $a,b,c,d$ etc. added to some multiple of $N$ shall
refer to this class.  
On the other hand,
all the remaining numbers less than $N$, whose total is also
$\frac{1}{2}(N-1)$, and which can be called {\em non-residues}, have
been given as the latter class, and designated by the Greek letters
\[
\alpha,\quad \beta,\quad \gamma,\quad \delta\quad \textrm{etc.}
\]
Concerning these two types of numbers, we have already noted above that the
product of two residues or members of the first class again falls in this class,
for example $ab,ac,bc$ etc., reducing them by division to be less than $N$,
and the product of a residue and a non-residue will appear in the latter
class of non-residues, and finally the product of two members of the non-residues
will again be a residue. With this noted, we prepare our demonstration so that
we would find a great absurdity to follow if no formula $\lambda xx+\mu yy+
\nu zz$ could be given that is divisible by $N$. The demonstration will
proceed in the following way.

I. Since all squares are equal to some residue $a$ or $b$ or $c$ added
to a particular multiple of the number $N$,
for the formula $\lambda xx+\mu yy+\nu zz$ were to be divisible by the number $N$, because $xx=\zeta N+a$,
$yy=\eta N+b$ and $zz=\vartheta N+c$,
it is certainly equivalent that
the formula $\lambda a+\mu b+\nu c$ be divisible by $N$.
If our theorem were false, there should be no way to make the formula
$\lambda a+\mu b+\nu c$ divisible by $N$.

II. Then since no formula of this type could be given that is divisible by $N$,
still less could it be $=0$, and thus this equation $\lambda a=-\mu b-\nu c$
will be impossible, and equally such an equation
\[
\lambda a=(\zeta N-\mu)b+(\eta N-\nu)c.
\]
Truly, because $\lambda,\mu$ and $\nu$ are prime to $N$, the coefficients
$\zeta$ and $\eta$ can always be chosen so that the formulas $\zeta N-\mu$
and $\eta N-\nu$ become divisible by $\lambda$. Let us therefore put
\[
\zeta N-\mu=\lambda m \quad \textrm{and} \quad \eta N-\nu=\lambda n
\]
and the equation
\[
a=mb+nc
\]
will be impossible too.

III. Therefore, since the formula $mb+nc$ cannot be equal to $a$, and hence
cannot appear in the class of residues (for admitting the contrary will negate
our theorem), 
it necessarily appears in the other class of non-residues;
there at once (because $c$ can denote unity)
$mb+n$ will occur, and hence all the formulas
\[
ma+n, \quad mb+n, \quad mc+n, \quad md+n \quad \textrm{etc.};
\]
since all these numbers are mutually distinct and there are $\frac{1}{2}(N-1)$
in total, they completely exhaust the class of non-residues,
of course dividing them by $N$ to make them below $N$.

IV. Indeed the products of all these numbers with any number of the first
class, such as $d$, will also occur in this class, which will therefore be
\[
mad+nd, \quad mbd+nd, \quad mcd+nd \quad \textrm{etc.}
\]
Indeed, the products $ab,bd,cd$ etc. fall in the first class and
appear among the
numbers $a,b,c,d$ etc.; and thus all the formulas\footnote{Translator: Since
$x \mapsto xd$ is a permutation of the residues, $mx+nd \mapsto mxd+nd$ is a 
permutation of the non-residues.}
\[
ma+nd, \quad mb+nd, \quad mc+nd \quad \textrm{etc.},
\]
will occur in the latter class among the non-residues,
which each exceed the preceding by
the quantity $n(d-1)$.\footnote{Translator: That is, the terms $ma+nd,mb+nd,mc+nd$ etc.
each
exceed the corresponding terms $ma+n,mb+n,mc+n$ etc. by $n(d-1)$.}
For the sake of brevity let
us put this difference $=\omega$, 
which will be prime to this divisor $N$ as long as $d$ is not assumed to be
unity, for $d-1$ is $<N$ and too the number $n$ is prime to $N$.

V. Thus if a number $\alpha$ is contained in the class of non-residues,
then simultaneously $\alpha+\omega$ will also occur, and for the same
reason this number again incremented by $\omega$, namely $\alpha+2\omega$,
and for the same reason the numbers $\alpha+3\omega, \alpha+4\omega$ etc.
also occur here. Therefore all the terms
of this arithmetic progression
\[
\alpha, \quad \alpha+\omega,\quad \alpha+2\omega,\quad \alpha+3\omega\quad \textrm{etc.},
\]
divided of course by $N$ to remain below $N$, will occur among the non-residues.

VI. Since the difference of this progression is $\omega$, namely a number 
prime to $N$, in this progression there will not only be a term divisible by $N$,
but will even yield all the numbers $1,2,3,4$ etc. with no exceptions when
divided by $N$. Thus it follows from the contrary hypothesis that
all the numbers whatsoever $1,2,3,4$ etc. occur in the class of non-residues;
since this is absurd, this contrary hypothesis is certainly false. Namely
it is false that no numbers of the form
\[
\lambda xx+\mu yy+\nu zz,
\]
can be given which are divisible by $N$. Therefore indeed such numbers
can be given;
and this is that
which we set out to show.

\begin{center}
{\Large Corollary 1}
\end{center}

Not only will always be possible to find
three squares $xx,yy$ and $zz$ of this type,
but also one of them, say $zz$, is our choice, providing it is not divisible
by $N$. Thus if $f$ denotes any number we please which is not divisible by $N$,
two squares $xx$ and $yy$ can always be assigned so that the formula
\[
\lambda xx+\mu yy+\nu ff
\]
becomes divisible by $N$. For demonstrating this, 
if $z$ is any umber, a number $v$ can always be given so that the product
$vz$ divided by $N$ leaves the given residue $f$. For let
$vz=\vartheta N+f$ and our formula multiplied by $vv$,
which certainly will still be divisible by $N$, would become
\[
\lambda vvxx+\mu vvyy+\nu(\vartheta\vartheta NN+2\vartheta Nf+ff),
\]
where, because the term $\vartheta \vartheta NN+2\vartheta Nf$ is itself divisible by $N$, 
the remaining form
\[
\lambda vvxx+\mu vvyy+\nu ff
\]
will be divisible by $N$.

\begin{center}
{\Large Corollary 2}
\end{center}

Whatever the numbers $\lambda, \mu, \nu$ are, unity or another number number can
always be assumed for one of them.
For since by multiplying by $\vartheta$ the formula
\[
\vartheta \lambda xx+\vartheta \mu yy+\vartheta \nu zz
\]
will admit division by $N$, in place of $\vartheta$ any such number can be chosen
so that the product $\vartheta \lambda$ leaves unity when divided
by $N$; for then the formula
\[
xx+\vartheta \mu yy+\vartheta \nu zz
\]
will even now be divisible by $N$. Furthermore, here we can write
the residues arising from division by $N$ in place of $\vartheta \mu$ and 
$\vartheta \nu$,
in which way we arrive at the same formula
which the Celebrated Lagrange considered.

\begin{center}
{\Large Scholion}
\end{center}

Thus behold that we have completed an unconditional demonstration of this
most notable theorem for all numbers, that all numbers whatsoever are the
sum of four or fewer squares, which indeed Fermat previously professed
to have found, but which has perished most sadly through the passage of the years.
There can certainly be no doubt that the demonstration of Fermat was much
simpler and more general than the one which now at last comes to light.
What makes it likely that his demonstration followed different lines is that
he demonstrated from the same source that all numbers are the sum
of three triangular numbers or fewer, then too
sums of five pentagonal numbers or fewer, and too sums of
six hexagonal numbers and so on, and this generality is altogether
missing from our conclusion.
And we are even still ignorant of a demonstration that
any number is a sum of three triangular numbers or fewer.
In the meanwhile however it is appropriate to observe about
this theorem that it is only true for integral numbers, unlike the 
other which we have demonstrated holds even for fractional
numbers;\footnote{Translator: Euler showed in Theorem
20, \S 97 of E242, {\em Demonstratio theorematis Fermatiani omnem numerum sive integrum sive fractum esse summam quatuor pauciorumve quadratorum}, that every rational number
is a sum of four squares of rational numbers.}
for all the fractions
$\frac{1}{2},\frac{3}{2},\frac{5}{2},\frac{7}{2},\frac{9}{2}$ etc.
do not allow themselves to be resolved into three triangular numbers in any way, that is one cannot find any rational values in place of $x,y,z$ such that
\[
\frac{1}{2}=\frac{xx+x}{2}+\frac{yy+y}{2}+\frac{zz+z}{2};
\]
whence, which seems rather remarkable, the equation
\[
1=xx+x+yy+y+zz+z
\]
is impossible for whatever fractional numbers are taken for $x,y,z$.

\end{document}